\documentclass{article}

\begin{document}

\begin{center}

\large\textbf{A Fourier approach to pizza inequity}

\end{center}

\bigskip

\bigskip

\begin{center}

\textbf{ABSTRACT}

\end{center}

\medskip

\noindent Let $n$ be an odd number greater than 1. We slice a circular pizza into
$2n$ slices, making cuts from a noncentral interior point of the circle.We estimate
the difference between the total area of the even numbered slices and the total
area of the odd numbered slices.

\bigskip

\bigskip

\bigskip

\bigskip

\bigskip

\noindent David Gluck \newline
Department of Mathematics\newline
Wayne State University\newline
Detroit, Michigan, 48202\newline
d.gluck@me.com

\bigskip

\bigskip

\bigskip

\bigskip

\bigskip

\bigskip


\bigskip

\bigskip

\bigskip

\bigskip


\newpage

\noindent \textbf{0. Introduction}

\medskip

To cut a circular pizza into $2n$ slices, one normally takes $2n$ rays
which begin at the center $C$ of the circle, such that the angle between
any two adjacent rays is $2\pi/2n$. The slices are the regions bounded
by two adjacent rays and the circle. We number the slices counterclockwise,
and so we have $n$ odd numbered slices and $n$ even numbered slices.
The so-called pizza theorem asserts, somewhat surprisingly, that even if
the $2n$ rays begin at the same noncentral interior point $P$ of the circle,
and $2n$ is divisible by 4 and greater than 4, then the total area of the even
numbered slices equals the total area of the odd numbered slices. According to the 
Wikipedia pizza theorem article, this fact originated as a problem in [3], and the
first published proof was given in [2].

Here we consider the case that $2n$ is not a multiple of 4. In this case, the pizza
theorem does not hold; see the interesting and intricate paper  [1]. We show here, however,
that the theorem is often almost true. We suppose without loss of generality that
the circle has radius 1, and we let $a$ be the distance from $P$ to $C$. Thus
$0<a<1.$ Let $\alpha$ denote the angle between the first ray and the line joining
$C$ to $P$. Let $g(\alpha,a,n)$ denote the total area of the even numbered slices
minus the total area of the odd numbered slices. We show in Corollary 1.5 that
$$|g(\alpha,a,n)|<a^n/2(1-a^2)(1-a^{2n}),$$
so  the ``pizza inequity'' is on the order of $a^n$ when $a$ is bounded away from 1.
This inequity is too small, one imagines, to be noticed by hungry pizza eaters.

Our key idea is to express an appropriate step function on $[0,2\pi]$ as a Fourier
series. This leads to Theorem 1.2, in which $g(\alpha,a,n)$ is given as an infinite
linear combination of functions $\sin m\alpha$, where $m$ ranges over the positive
odd multiples of $n$. The coefficient of $\sin m\alpha$ is a power series in $a$.
The (numerical) coefficients in this power series all have the same sign. We evaluate
these power series coefficients explicitly in terms of binomial coefficients. We also
find an exact expression for the maximum value of $|g(\alpha,a,n)|$ when $a$ and 
$n$ are fixed.

\medskip

\noindent \textbf{1. Results}

\medskip

\noindent \textit{Notational conventions.} In this paper, $n$ denotes a fixed odd integer
greater than 1, and $a$ denotes a real number in the interval (0,1). The notation
$\sum_m$ means that $m$ ranges over all positive odd multiples of $n$. For a real number
$x$, $e(x)$ denotes $e^{ix}$. If $r$ is a positive integer or $r=1/2$, and $s$ is an integer
we define $C(r,s)$ to be the binomial coefficient ``$r$ choose $s$''. As usual, this is
defined to be zero if $s<0$ or if $r$  is an integer and $s>r$. If $r$ and $s$ are
integers, then $\delta_{r,s}$ is the Kronecker delta.

\medskip

Consider the circle $(x-a)^2+y^2=1$, with center ($a$,0) and radius 1. Thus ($a$,0) and (0,0)
respectively play the roles of $C$ and $P$ in the introduction. Fix an initial angle $\alpha \in
[0, 2\pi)$ and consider the $2n$ rays starting from the origin and making angles $\alpha,
\alpha+\pi/n, \ldots, \alpha +(2n-1)\pi/n$ with the positive $x$-axis. These rays divide the 
pizza into $2n$ slices. We number them counterclockwise, so the first slice is bounded by
the $\alpha$ ray, the $\alpha +\pi/n$ ray, and the circle. The second slice is bounded by the
$\alpha+\pi/n$ ray, the $\alpha+2\pi/n$ ray, and so on. Let $A$ denote the union of the $n$
even numbered intervals $[\alpha+\pi/n,\alpha+2\pi/n], \ldots$ and let $B$ denote the union of the $n$
odd numbered intervals $[\alpha, \alpha+\pi/n], \ldots$.

In polar coordinates, the circle has equation $r(\theta)=a\cos\theta+\sqrt{1-a^2\sin^2 \theta}$;
taking the negative sign in the quadratic formula would cause $r$ to be negative. We define the
pizza inequity function $g(\alpha,a,n)$ to be the total area of the even numbered slices minus
the total area of the odd numbered slices. Hence $$2g(\alpha,a,n)=\int_A r(\theta)^2 d\theta
-\int_B r(\theta)^2 d\theta.$$ Squaring our expression for $r(\theta)$ yields $2g(\alpha,a,n)=$
$$\int_A1+a^2\cos2\theta+2a\cos\theta\sqrt{1-a^2\sin^2\theta} d\theta$$
\newline
$$- \int_B1+a^2\cos2\theta+2a\cos\theta\sqrt{1-a^2\sin^2\theta} d\theta.$$

We remark that if we replaced our circle $(x-a)^2+y^2=1$ by the circle with center $(aR,0)$
and radius $R$, then the function $r(\theta)$ would be multiplied by $R$ and so the integrands in the
last two integrals would be multiplied by $R^2.$ Hence the ratio of the pizza inequity to the
area of the circle would be unchanged. Thus we lose no generality by working with a circle
of radius 1. 

\medskip

\noindent\textbf{Lemma 1.1} Let $n, a, \alpha, A,$ and $B$ be as above. Let $s(\theta)$ denote
the step function with values +1 on $A$ and -1 on $B$; the values of $s(\theta)$ at the endpoints 
of the $2n$ intervals don't matter. Let $m$ be an integer. Then $\int_0^{2\pi} s(\theta)e(-m\theta) d\theta$
=0 unless $m$ is an odd multiple of $n$. If $m$ is an odd multiple of $n$, then
$$\int_0^{2\pi} s(\theta)e(-m\theta) d\theta=(-4n/mi)e(-m\alpha).$$

\medskip

\noindent\textit{Proof.} We may assume that $m\ne 0.$ Since $-e(-m\theta)/mi$ is an antiderivative of $e(-m\theta)$, we see that
 $\int_A e(-m\theta) d\theta$ equals 
 $$-(1/mi)[e(-m(\alpha +2\pi/n))+e(-m(\alpha +4\pi/n))+\ldots +e(-m\alpha)]$$
  $$+(1/mi)[e(-m(\alpha+\pi/n))+e(-m(\alpha+3\pi/n))+\ldots+e(-m(\alpha+(2n-1)\pi/n))].$$
 Let $\zeta=e(2\pi/n)$ and let $\eta=e(\pi/n)$. The foregoing becomes
 $$-(e(-m\alpha)/mi)[\zeta^{-m}+\zeta^{-2m}+\ldots+1]$$
  $$+(e(-m\alpha)/mi)[\eta^{-m}(1+\zeta^{-m}+\ldots+\zeta^{-(n-1)m})].$$
 If $n$ does not divide $m$, then the map $\zeta\mapsto\zeta^{-m}$ defines a nontrivial character
  of the multiplicative group generated by $\zeta$, and so both expressions in brackets vanish. If $2n$
  divides $m$, then both expressions in brackets equal $n$ and again the integral vanishes. Suppose
  next that $m$ is an odd multiple of $n$. Then $\eta^{-m}=-1$ and so $\int_A e(-m\theta) d\theta=
  (-2n/mi)e(-m\alpha).$
  
  A similar computation shows that $\int_B e(-m\theta) d\theta$ equals
  $$-(e(-m\alpha)/mi)[\eta^{-m}(1+\zeta^{-m}+\zeta^{-2m}+\ldots+\zeta^{-(n-1)m})]$$
   $$+(e(-m\alpha)/mi)[\zeta^{-m}+\zeta^{-2m}+\ldots+1].$$
  As before, we see that $\int_B e(-m\theta) d\theta$ vanishes unless $m$ is an odd multiple of
  $n$, in which case $\int_B e(-m\theta) d\theta=+(2n/mi)e(-m\alpha).$ Since
  $$\int_0^{2\pi} s(\theta)e(-m\theta) d\theta=\int_A e(-m\theta) d\theta-\int_B e(-m\theta) d\theta,$$
  the conclusion of the lemma follows.
  
  \medskip
  
  \noindent \textbf{Theorem 1.2} For $n>1$ odd, $0<a<1$, and $\alpha \in [0,2\pi)$, let
  $g(\alpha,a,n)$ be the pizza inequity function defined above. Let$f(\alpha,a,n)=g(\alpha,a,n)/a.$
  Then $$f(\alpha,a,n)=\sum_m (4n/\pi m)(P_m(a))\sin m\alpha,$$
  where $P_m(x)$ is the power series $\sum_{j=1}^{\infty} c_{2j}(m)x^{2j},$ with
  $$c_{2j}(m)=(-1)^jC(1/2,j)\int_0^{2\pi} \cos \theta \cos m\theta \sin^{2j} \theta d\theta.$$
  
  \medskip
  
  \noindent\textit{Proof.} The step function $s(\theta)$ of Theorem 1.1 has  Fourier series
  $\sum_{m=0}^{\infty} a_m\cos m\theta+\sum_{k=1}^{\infty} b_k\sin k\theta$. Since
  $\int_0^{2\pi} s(\theta)d\theta=0,$ we have $ a_0=0.$ For $m>0,$
  $$a_m=(1/\pi)\int_0^{2\pi} s(\theta)\cos m\theta d\theta=(1/2\pi)\int_0^{2\pi}s(\theta)(e(-m\theta)+
  e(m\theta)) d\theta.$$
  By Lemma 1.1, if $m$ is an odd multiple of $n$, then
  $$\int_0^{2\pi} s(\theta)(e(-m\theta)+e(m\theta)) d\theta=(-4n/mi)e(-m\alpha)+(4n/mi)e(m\alpha)$$
  $$=(4n/mi)(e(m\alpha)-e(-m\alpha))=(8n/m)(1/2i)(e(m\alpha)-e(-m\alpha))=(8n/m)\sin m\alpha.$$
  Hence $a_m=(4n/\pi m)\sin m\alpha.$ If $m$ is not an odd multiple of $n$, then $a_m=0$ by Lemma 1.1.We 
  won't bother to compute $b_k$, since we will soon see that its value doesn't matter.
  
  Now 
  $$2g(\alpha,a,n)=\int_0^{2\pi}s(\theta)\left[1+a^2\cos 2\theta+2a\cos \theta \sqrt{1-a^2\sin ^2\theta}\right] d\theta.$$
  The functions 1 and $\cos 2\theta$ are orthogonal on $[0,2\pi]$ to $\cos m\theta$, whenever $m$ is an odd multiple of $n$.
  The functions 1 and $\cos 2\theta$ are also orthogonal to any function $\sin k\theta$. Hence
  $$g(\alpha,a,n)=a\int_0^{2\pi}s(\theta)\cos \theta\sqrt{1-a^2\sin^2 \theta} d\theta.$$
  If $k$ is an integer, then
  $$\int_0^{2\pi}\sin k\theta \cos \theta \sqrt{1-a^2\sin ^2 \theta}d\theta=\int_{-\pi}^{\pi}\sin k\theta \cos \theta
  \sqrt{1-a^2\sin^2\theta} d\theta=0,$$
  since the last integrand is an odd function of $\theta$. Hence
  $$f(\alpha,a,n)=\sum_m(4n/\pi m)\left[ \int_0^{2\pi}\cos \theta \cos m\theta \sqrt{1-a^2\sin^2 \theta} d\theta \right]\sin m\alpha.$$
  
  The binomial series
  $$(1-x)^{1/2}=\sum_{j=0}^{\infty}C(1/2, j)(-x)^j=\sum_{j=0}^{\infty}(-1)^jC(1/2,j)x^j$$
  converges for $|x|<1.$ Note that all coefficients in this series, except the constant term, are negative.
  The series $\sum_{j=0}^{\infty} (-1)^jC(1/2,j)(a^2\sin^2\theta)^j$ converges uniformly in $\theta$ to
  $\sqrt{1-a^2\sin^2\theta}.$ Thus
  $$\int_0^{2\pi} \cos \theta \cos m\theta \sqrt{1-a^2\sin^2\theta} d\theta=\sum_{j=0}^{\infty} (-1)^jC(1/2,j)
  \left [\int_0^{2\pi} \cos \theta \cos m\theta \sin^{2j}\theta d\theta\right ]a^{2j}.$$
  Since $m\geq n>1$, we have $\int_0^{2\pi} \cos \theta \cos m\theta d\theta=0,$ so we can start the
  summation at $j=1.$ The conclusion of the theorem follows.
  
  \medskip
  
  \noindent\textbf{Proposition 1.3} Let $P_m(x)=\sum_{j=1}^{\infty} c_{2j}(m)x^{2j}$ as in Theorem 1.2. Then
   $$c_{2j}(m)=(-1)^{(m+1)/2} (\pi /2^{2j}) |C(1/2,j)|[C(2j, (2j-m+1)/2)-C(2j, (2j-m-1)/2)].$$
   All nonzero coefficients $c_{2j} (m)$ have the same sign, namely $(-1)^{(m+1)/2}$. The leading
   coefficient of $P_m(x)$ is
   $$c_{m-1}(m)=(-1)^{(m+1)/2}(\pi/2^{m-1})|C(1/2,(m-1)/2)|$$
   For all $m$ and $j$ we have
   $$|c_{2j}(m)|\leq |c_2(3)|=\pi/8.$$
   
   \medskip
   
   \noindent \textit{Proof.} Since $(-1)^jC(1/2,j)<0$ for $j\geq 1$, Theorem 1.2 implies that the first
   assertion amounts to saying that $\int_0^{2\pi} \cos \theta \cos m\theta \sin^{2j} \theta d\theta$ equals
   $$(-1)^{(m-1)/2}(\pi/2^{2j})
   [C(2j, (2j-m+1)/2)-C(2j, (2j-m-1)/2)].$$ 
   Now $\int_0^{2\pi} \cos \theta \cos m\theta \sin^{2j}\theta d\theta$ equals
   $$i^{2j}2^{-(2j+2)}\int_0^{2\pi} (e(\theta)+e(-\theta))(e(m\theta)+e(-m\theta))(e(\theta)-e(-\theta)^{2j} d\theta.$$
   Since $i^{2j}=(-1)^j$, this equals $$2^{-(2j+2)}\sum_{k=0}^{2j}(-1)^{j+k}C(2j,k)I_k,$$ with
   $$I_k=\int_0^{2\pi}[e((m+1)\theta)+e((-m-1)\theta)+e((m-1)\theta)+e((1-m)\theta)]
   e((2k-2j)\theta) d\theta$$
   By orthogonality of exponential functions, 
   $$I_k=2\pi [\delta_{2j-2k,m+1}+\delta_{2j-2k,m-1}+\delta_{2j-2k,-m-1}+\delta_{2j-2k, 1-m}].$$ 
   If $2j-2k=\pm (m+1)$, then $k=(2j\mp(m+1))/2$ and
   $$(-1)^{j+k}C(2j,k)=(-1)^{(m+1)/2}C(2j,(2j-m-1)/2).$$
   If $2j-2k=\pm (m-1)$, then $k=(2j\mp (m-1))/2$ and
   $$(-1)^{j+k}C(2j,k)=(-1)^{(m-1)/2}C(2j, (2j-m+1)/2).$$ It follows that
   $$\sum_{k=0}^{2j} (-1)^{j+k}C(2j, k)I_k=4\pi(-1)^{(m-1)/2}[C(2j,(2j-m+1)/2)-C(2j, (2j-m-1)/2)]$$
   which, when multiplied by $2^{-(2j+2)}$, gives the desired value for our integral.
   
   Since $(2j-m-1)/2<(2j-m+1)/2<j$, the second assertion of the proposition follows immediately
   from the first. Since $(2j-m+1)/2<0$ if $2j<m-1,$ it follows that $c_{m-1}(m)$ is the leading
   coefficient of $P_m(x).$ Since $C(2j,0)-0=1,$ the third assertion of the proposition follows.
   
   To prove the final assertion, first suppose that $j \geq 3.$ For a positive integer $t$, let
   $M(t)$ denote the maximum ratio $C(t,u)/2^t$ for $u=0,1,\ldots,t.$ If $0\leq v\leq t+1,$
   then $C(t+1,v)=C(t,v-1)+C(t,v)$. It follows easily that $C(t+1,v)/2^{t+1} \leq M(t),$
   and so $M(t+1) \leq M(t).$ In particular, if $t \geq6, $ then $M(t)\leq M(6)=5/16.$
   Hence the factor 
   $$2^{-2j}[C(2j,(2j-m+1)/2)-C(2j,(2j-m-1)/2]$$
   of $c_{2j}(m)$ is at most 5/16, and so
   $$|c_{2j}(m)| \leq(5\pi/16)|C(1/2,3)|=(5\pi/16)(1/16)<\pi/8.$$
   
   If $j=2$, then $c_{2j}=0$ or $m$ is 3 or 5. One computes that $|c_4(5)|=\pi/128$
   and $|c_4(3)|=3\pi/128$, both less than $\pi/8.$ If $j=1,$ then $c_{2j}=0$ or $m=3.$
   One computes that $|c_2(3)|=\pi/8,$ as desired. 
   
   \medskip
   
   \noindent\textit{Remark.} Our bound for $|c_{2j}(m)|$ could be improved considerably if
   one assumes, say, that $n\geq 5.$ We leave this to the interested reader.
   
   \medskip
   
   \noindent\textbf{Lemma 1.4} For $\alpha,a,$ and $n$ as above, let $f_a(\alpha)=f(\alpha,a,n).$
   Then $f_a$ is an odd function of $\alpha$ and $f_a$ is periodic with period $2\pi/n.$
   
   \medskip
   
   \noindent\textit{Proof.} This is immediate from Theorem 1.2.
   
   \medskip
   
   \noindent\textbf{Corollary 1.5} With notation as above, let $M_a=\sum_m (4n/\pi m)|P_m(a)|.$
   Then $f_a(\pi/2n)=(-1)^{(m+1)/2}M_a$ and $f_a(-\pi/2n)=(-1)^{(m-1)/2}M_a.$ The maximum and 
   minimum values of $f_a$ are $M_a$ and $-M_a,$ respectively. We have
   $$0<M_a<a^{n-1}/2(1-a^2)(1-a^{2n}).$$
   Consequently, for all $\alpha,$ the true pizza inequity function $g(\alpha,a,n)$ satisfies
   $$|g(\alpha,a,n)|<a^n/2(1-a^2)(1-a^{2n}).$$
   
   \medskip
   
   \noindent\textit{Proof.} Note that $|P_m(a)|=\pm P_m(a)$ by Proposition 1.3. By Theorem 1.2,
   $$f_a(\pi/2n)=(4/\pi)P_n(a)-(4/3\pi)P_{3n}(a)+(4/5\pi)P_{5n}(a)-\ldots$$
   First suppose that $n\equiv 1$ (mod 4). By Proposition 1.3, all nonzero coefficients in $P_n(x)$
   are negative, all nonzero coefficients in $P_{3n}(x)$ are positive, all nonzero coefficients in
   $P_{5n}(x)$ are negative, and so on. It follows that $f_a(\pi/2n)=-M_a$. If $n\equiv 3$ (mod 4),
   similar reasoning shows that $f_a(\pi/2n) =+M_a.$ By Lemma 1.4, $f_a(-\pi/2n)=-f_a(\pi/2n).$
   By the triangle inequality, $|f_a(\alpha)|\leq M_a$ for all $\alpha$. Hence we have found the 
   maximum and minimum values of $f_a.$
   
   We now estimate $M_a.$ By Proposition 1.3, we have $|c_{2j}(m)| \leq \pi/8$ for all $j$ and $m.$
   Since $c_{2j}=0$ for $2j<m-1,$ comparison with a geometric series shows that 
   $$|P_m(a)| \leq (\pi/8)a^{m-1}/(1-a^2)$$
   for all $m.$ Hence
   $$M_a \leq \sum_m(4n/\pi m)(\pi/8)a^{m-1}/(1-a^2)<\sum_ma^{m-1}/2(1-a^2).$$
   But the last series is geometric with ratio $a^{2n}$ and first term $a^{n-1}/2(1-a^2).$
   Hence $M_a<a^{n-1}/2(1-a^2)(1-a^{2n}),$ as desired. Since $g(\alpha,a,n)=af(\alpha,a,n),$
   the final assertion of the corollary follows.
   
   \newpage
   
   \begin{center}{\textbf{References}}\end{center}
   
   \medskip
   
   \noindent [1] P. Deiermann and R. Mabry, Of cheese and crust: a proof of the pizza
   conjecture and other tasty results, Amer. Math. Monthly (2009), 423-438.
   
   \smallskip
   
   \noindent [2] M. Goldberg, Divisors of a circle, Math. Mag. 41(1968), 46.
   
   \smallskip
   
   \noindent [3] L. J. Upton, Problem 660, Math. Mag. 40(1967), 163.

\end{document}